\documentclass[Preprint,12pt]{elsarticle}

\usepackage{color,amsmath,latexsym,amssymb}
\usepackage{amsfonts}

\newcommand{\B}{\mathbb{B}}
\newcommand{\E}{\mathbb{E}}
\newcommand{\G}{\mathbb{G}}
\newcommand{\W}{\mathbb{P}}

\newcommand{\R}{\mathbb{R}}
\newcommand{\N}{\mathbb{N}}

\newcommand{\cA}{\mathcal{A}}
\newcommand{\cF}{\mathcal{F}}

\newcommand{\cP}{\mathcal{P}}

\newcommand{\fF}{\mathfrak{F}}

\def \one {\mathchoice {\hbox{1\kern -0.27em l}}{\hbox{1\kern -0.27em l}}
                     {\small{1\kern -0.27em l}}{\small{1\kern -0.27em l}}}
                     
\newtheorem{thm}{Theorem}                     
\newtheorem{lem}[thm]{Lemma}
\newtheorem{prop}[thm]{Proposition}
\newtheorem{cor}[thm]{Corollary}
\newproof{pf}{Proof}
\newdefinition{rmk}{Remark}
\newdefinition{exam}{Example}

\begin{document}

\begin{frontmatter}

\title{Asymptotic Properties of Monte Carlo Methods\\ 
in Elliptic PDE-Constrained Optimization\\ 
under Uncertainty} 

\author[1]{W.~R\"omisch\corref{cor1}%
           }
\ead{romisch@math.hu-berlin.de}           

\author[2]{T.~M. Surowiec}
\ead{surowiec@mathematik.uni-marburg.de}

\affiliation[1]{organization={Institute of Mathematics, 
                Humboldt University of Berlin},\\
                addressline={Unter den Linden 6},
                postcode={10099}, 
                city={Berlin},
                country={Germany}} 
                
\affiliation[2]{organization={FB12 Mathematik und
                Informatik, Philipps-Universit\"at Marburg},
                addressline={Hans-Meerwein-Stra{\ss}e 6},
                postcode={35032}, 
                city={Marburg},
                country={Germany}}
  		 
\cortext[cor1]{Corresponding author}

\begin{abstract}
Monte Carlo approximations for random linear elliptic 
PDE constrained optimization problems are studied. We 
use empirical process theory to obtain best possible 
mean convergence rates $O(n^{-\frac12})$ for optimal 
values and solutions, and a central limit theorem for 
optimal values. The latter allows to determine 
asymptotically consistent confidence intervals by 
using resampling techniques.
\end{abstract}

\begin{keyword}
random elliptic PDE \sep stochastic optimization \sep 
Monte Carlo \sep central limit theorem \sep resampling
\MSC 49J20 \sep 49J55 \sep 60F17 \sep 65C05 \sep
90C15 \sep 35R60
\end{keyword}

\end{frontmatter}

%\maketitle

\section{Introduction}
\label{sec:intro}
PDE-constrained optimization under uncertainty is a rapidly 
growing field with a number of recent contributions in 
theory \cite{ChQR13,KoSu16,KoSu18,KoSu20}, 
numerical and computational methods
 \cite{GaSU21,GeWo20,KHRW13,vBVa19}, and applications 
\cite{ChNK20,ChHG21,CHPRS11,ScSS11}. Nevertheless, 
a number of open questions remain unanswered, even for 
the ideal setting including a strongly convex objective 
function, a closed, bounded and convex feasible set, and 
a linear elliptic PDE with random inputs.

Broadly speaking, the numerical solution methods available
for such problems derive in part from the standard 
paradigms found in the classical stochastic programming
literature: first-order stochastic 
approximation/stochastic gradient (SG) approaches,
\cite{GeWo20}, versus methods that rely on sampling from 
the underlying probability measure \cite{vBVa19}. The 
latter approaches are not optimization algorithms per 
se, rather, they use approximations of the expectations 
by replacing the underlying probability distribution with 
a discrete measure. This may be obtained either from
available data or using Monte Carlo (MC) type methods. 
The main advantage of the latter is that we may turn to 
the wide array of powerful, function-space-based 
numerical methods for PDE-constrained optimization in 
a deterministic setting, see e.g., \cite{HPUU09}. 

Nevertheless, as with SG-based methods, the optimal 
values and solutions of MC-type approximations must also 
be understood as realizations of a rather complex random 
process.  In this context, it is helpful to think of them 
as mappings from some space of probability measures into 
the reals (optimal values) and decision space (solutions).  
In optimization under uncertainty, stability usually 
refers to the continuity properties of these mappings 
with respect to changes in the underlying measure. As 
the underlying parameter space is not a normed linear
space, proving continuity and asymptotic
properties as $n\to\infty$ can be a delicate matter.
Such statements require techniques not typically 
employed in PDE-constrained optimization, e.g., empirical
process theory or the method of probability metrics.
Of course, if we can obtain computable quantitative bounds 
in $n$, i.e., convergence rates, then such stability
statements can provide us with a priori information 
regarding the necessary sample size for an MC-based
numerical solution method. These can in turn be linked 
to the PDE discretization error for a comprehensive a 
priori error estimate.

In a recent paper, we provided a number of qualitative 
and quantitative stability statements for 
infinite-dimensional stochastic optimization problems 
using the method of probability metrics \cite{HoRS20}. 
However, for the PDE-constrained optimization problem
provided in \cite[Sec. 7]{HoRS20}, a major open question
remained: Can we derive a reasonable rate of convergence 
for the minimal information metric that supports our
numerical results?  Taking this question as a starting 
point, we seek to answer this by using deep results from 
empirical process theory as detailed in
\cite{Gine96,vdVWe96}. The idea to use empirical 
process theory has been employed in the stochastic
programming literature before, cf. some results in
\cite[Chapters 6--8]{RuSh03}. However, it has not been 
used in situations where the decision spaces are infinite
dimensional, which presents an additional challenge.
Another approach based on large deviation-type
results is employed in the recent preprint \cite{Milz21}
in which the author obtained results that are in parts
similar those in the present paper. However, the results
in \cite{Milz21} cannot be used to derive confidence
intervals for the optimal values without further 
assumptions on the integrands.

Convergence statements based on the underlying smoothness 
of the uncertainty in the forward-problem have been
previously considered in \cite{KuSc13}. However, 
these results do not in fact subsume those found in the
present article. There are a number of significant 
differences. First, the authors in \cite{KuSc13} do not 
place additional constraints on the control variables.
Therefore, the first-order optimality conditions are 
coupled systems of PDEs. Second, the authors use this 
fact to transform the question of convergence rates of 
the optimization problem into convergence rates for the 
PDE system. Following this, they allow the controls to 
also depend on the uncertainty, but they do not include 
an additional non-anticipativity constraint of the type: 
$u(\sigma)= \E_{\sigma}[u]$ for all parameters $\sigma$.
Therefore, the rates no longer apply to the original
optimization problem. As a consequence, our results 
appear to be the first of their kind for PDE constrained
optimization under uncertainty.

The rest of the article is organized as follows. In 
Section \ref{sec:problem} we introduce the PDE constrained
optimization problem \eqref{sopde} with random parameters 
studied in this article. In addition, we consider a
suitable distance measure for probability distributions.
After discussing its basic properties, we recall some 
results from \cite{HoRS20} on continuity properties of
infima and solutions to the stochastic optimization with 
respect to such probability metrics. Finally, we prove 
a new result on Lipschitz continuity for solutions.
Section \ref{sec:mc} contains the main results on
convergence rates for infima and solutions of Monte 
Carlo approximations and a central limit theorem for
infima. All results are consequences of empirical process
theory. In Section \ref{sec:subsampl} we shortly describe
how subsampling methods can be used to complement the
central limit result by deriving confidence intervals 
for the optimal values. We close the paper by discussing 
the limitations and possible extensions of our results.

\section{PDE constrained optimization under uncertainty}
\label{sec:problem}

We start by introducing several function spaces used
throughout our study. Let $(\Omega,\cA,\mu)$ be a 
measure space with $\sigma$-finite measure and $Y$ be 
a linear space endowed with a norm $\|\cdot\|$. 
By $L_p(\Omega,\cA,\mu;Y)$, $1\le p\le\infty$,
we denote the linear space of strongly measurable 
functions $y:\Omega\to Y$ such that the integral 
$\int_{\Omega}\|y(\omega)\|^{p}d\mu$ is finite. As 
usual we equip this space with the norm
\[
\|y\|_p=\left(\int_{\Omega}\|y(\omega)\|^{p}d\mu
\right)^{\frac{1}{p}}
\]
for $p<\infty$ with the usual modification for 
$p=\infty$. In case $V=\R$ we will omit $V$ and 
write  $L_p(\Omega,\cA,\mu)$. If $\Omega$ is a subset 
of some Euclidean space $\R^r$ and $\mu=\lambda$ 
the Lebesgue measure, we will shortly write 
$L_{p}(\Omega)$. If $\Omega$ is a subset of a metric 
space and $\cA$ the Borel $\sigma$-fileld, we will omit 
$\cA$. For any set $Y$ we denote by $\ell^{\infty}(Y)$ 
the linear space of bounded real-valued functions $g$ 
defined on $Y$ endowed with the norm 
$\|g\|=\sup_{y\in Y}|g(y)|$.
 
Now, let  $D\subset\R^m$ be an open bounded domain 
with Lipschitz boundary and $V=H_0^1(D)$ the usual 
Sobolev space of (equivalence classes of) functions in
$L_{2}(D)$ that admit square integrable weak derivatives. 
We endow this space with the  inner product 
$(u,v)_{V}=\int_{D}\nabla u \cdot \nabla v\, dx$ 
and norm $\|u\|=\sqrt{(u,u)}$. The topological dual 
is denoted by $V^{\star}=H^{-1}(D)$ with the usual 
operator norm $\|\cdot\|_{\star}$. In addition, we 
consider the Hilbert space $H=L^{2}(D)$ with the 
inner product $(g,h)_H=\int_{D}g(x)h(x)\,dx$. 
The dual pairing for $V,V^{\star}$ is denoted by 
$\langle\cdot,\cdot\rangle$.\\
Let $\Xi$ be a metric space and $\cP(\Xi)$ be the set of 
all Borel probability measures on $\Xi$. Fix $P\in\cP(\Xi)$.
For the parametric PDE, we first define the bilinear 
form $a(\cdot,\cdot;\xi):V\times V\to\R$ 
\begin{equation}\label{bilinform}
a(u,v;\xi)=\int_{D}\sum_{i,j=1}^{m}b_{ij}(x,\xi)
\frac{\partial u(x)}{\partial x_{i}}
\frac{\partial v(x)}{\partial x_{j}}\,dx
\end{equation}
for each $\xi\in\Xi$. Here, we impose the condition that 
the functions $b_{ij}:D\times\Xi\to\R$ are  measurable 
on $D\times\Xi$ and there exist $L>\gamma>0$ such that
\begin{equation}\label{unifellipt}
\gamma\sum_{i=1}^{m}y_{i}^{2}\le
\sum_{i,j=1}^{m}b_{ij}(x,\xi)y_iy_j\le L\sum_{i=1}^{m}
y_{i}^{2}\quad(\forall y\in\R^{n})
\end{equation}
for a.e. $x\in D$ and any $\xi\in\Xi$. This implies that
each $b_{ij}$ is essentially bounded on $D\times\Xi$ in 
both arguments with respect to the product measure
$\lambda\times P$.

We consider the stochastic optimization problem: 
Minimize the functional
\begin{eqnarray}\nonumber
{\mathcal{J}}(u,z)&:=&\frac{1}{2}\int_{\Xi}\int_{D} 
|u(x,\xi) - \widetilde{u}(x)|^2\,dx\,dP(\xi) 
+ \frac{\alpha}{2} \int_{D} |z(x)|^2 \,dx \\
&=& \frac{1}{2} \E_{P}[\|u(\cdot) - \widetilde{u}\|^2_{H}] 
+ \frac{\alpha}{2} \| z \|^2_{H}
\label{object}
\end{eqnarray}
subject to $(u,z)\in L_2(\Xi,P;V)\times Z_{\rm ad}$, 
where  $\alpha>0$, $\widetilde{u}\in H$, $Z_{\rm ad}$
denotes a closed convex bounded subset of $H$ and 
$u(\cdot)$ solves the random elliptic PDE
\begin{equation}\label{rellpde}
a(u(\xi),v;\xi)=\int_{D}(z(x)+g(x,\xi))v(x)\,dx 
\end{equation}
for $P$-a.e. $\xi\in\Xi$ and all test functions $v\in V$, 
where $g:D\times\Xi\to\R$ is measurable on $D\times\Xi$ 
and $g(\cdot,\xi)\in H$ for each $\xi\in\Xi$.

For $P$-a.e. $\xi\in\Xi$ we define the mapping 
$A(\xi):V\to V^{\star}$ by means of the Riesz 
representation theorem
\[
\langle A(\xi)u,v \rangle = a(u,v;\xi)\quad(u,v\in V).
\]
Consequently $A(\xi)$ is linear, uniformly positive 
definite (with $\gamma>0$) and uniformly bounded 
(with $L>0$) and the random elliptic PDE may be 
written in operator form
\[
A(\xi)u=z+g(\xi)\quad(P\mbox{-a.e. }\xi\in\Xi).
\]
In addition, the inverse mapping 
$A(\xi)^{-1}:V^{\star}\to V$ exists, and is linear, 
uniformly positive definite (with modulus $L^{-1}$) and
uniformly bounded (with constant $\gamma^{-1}$).
This allows us to rewrite the stochastic optimization 
problem in reduced form over $z\in Z_{\rm ad}$ as:
\begin{equation}\label{sopde}
\min\left\{F_{P}(z)=\int_{\Xi}f(z,\xi)\,dP(\xi): 
z\in Z_{\rm ad}\right\}
\end{equation}
with the integrand
\begin{eqnarray}\nonumber
f(z,\xi)&=&\frac12\big\|A(\xi)^{-1}(z+g(\xi))-
\widetilde{u}\big\|_{H}^{2}+\frac{\alpha}{2}\|z\|_{H}^{2}\\
&=&\frac{1}{2}\left\|A(\xi)^{-1}z -(\widetilde{u}-
A(\xi)^{-1}g(\xi))\right\|_{H}^{2}+\frac{\alpha}{2}
\|z\|_{H}^{2}
\label{sopde1}
\end{eqnarray}
for any $z\in H$ and $\xi\in\Xi$, where 
$g\in L_{2}(\Xi,P;H)$ and $A(\xi)^{-1}$ as defined 
earlier. For each $\xi\in\Xi$ the function 
$f(\cdot,\xi):H\to\R$ is convex and continuous.
For later use we denote the optimal value
of \eqref{sopde} by $v(P)$.\\
We will need a few properties of the function
$F_{P}:H\to\R$. They are collected in the following
result which is partly proved in \cite{HoRS20}.

\begin{prop}\label{p1}
For each $P\in\cP(\Xi)$ the functional $F_{P}$ is finite, 
continuous, and strongly convex on $H$ and, hence, weakly 
lower semicontinuous on the weakly compact set $Z_{\rm ad}$.  
Moreover, there exists an unique minimizer $z(P)\in 
Z_{\rm ad}$ of \eqref{sopde} and the objective function 
$F_{P}$ has quadratic growth around $z(P)$, i.e. we have
\begin{equation}\label{quadrgrowth}
\|z-z(P)\|_{H}^{2}\leq\frac{8}{\alpha}(F_{P}(z)-F_{P}(z(P)))
=\frac{8}{\alpha}(F_{P}(z)-v(P))\quad(\forall z\in 
Z_{\rm ad}).
\end{equation}
In addition, $F_{P}$ is G\^ateaux differentiable on $H$
with G\^ateaux derivative $F'_{P}(\cdot)$ and the estimate
\begin{equation}\label{Gdiff}
|F_{P}(z)-F_{P}(\tilde{z})|\le\sup_{t\in[0,1]}
\|F'_{P}(z+t(\tilde{z}-z))\|\,\|z-\tilde{z}\|
\end{equation}
holds for all $z,\tilde{z}\in H$.
\end{prop}
\begin{pf}
While the first part is proved in \cite{HoRS20}, it remains
to prove the G\^ateaux differentiability of $F_{P}$ and
the estimate \eqref{Gdiff}. For any $z,w\in H$ we
observe that the real function $h(t)=F_{P}(z+tw)$ is 
quadratic for $t\in\R$ and, hence, differentiable at
$t=0$. This means that $F_{P}$ is G\^ateaux 
differentiable at $z$. Now, we set $w=\tilde{z}-z$ for
some $\tilde{z}\in H$. Since a differentiable function 
$h:[0,1]\to\R$ satisfies the estimate
\[
|h(1)-h(0)|\le\sup_{t\in[0,1]}|h'(t)|,
\]
we obtain for $h(t)=F_{P}(z+t(\tilde{z}-z))$ the desired
estimate \eqref{Gdiff}.\hfill$\Box$
\end{pf}
Motivated by \eqref{sopde} and \eqref{quadrgrowth} we
consider the pseudo-metric
\begin{equation}\label{zetadist}
d_{\fF}(P,Q)=\sup_{f\in\fF}\left|\int_{\Xi}f(\xi)
\,dP(\xi) - \int_{\Xi}f(\xi)\,dQ(\xi)\right|
\end{equation}
on $\cP(\Xi)$ for studying quantitative stability of 
\eqref{sopde} with respect to perturbations of the 
underlying probability distribution $P$. Here, $\fF$ 
is a class of real-valued Borel measurable functions on 
$\Xi$. The notion of pseudo-metric means that all 
properties of metrics are satisfied except that 
$d_{\fF}(P,Q)=0$ does not imply $P=Q$ in general, 
unless $\fF$ is sufficiently rich. These are the 
typical properties required for probability metrics 
(see \cite{Rach91}). Such distances of probability 
measures were first introduced and studied in \cite{Zolo83}.
A number of important probability metrics are of the form 
$d_{\fF}$, for example, the bounded Lipschitz metric 
and the Fortet-Mourier metrics for which $\fF$ contains
(locally) Lipschitz functions. In both cases the
class $\fF$ is rich enough and, in addition, convergence 
with respect to $d_{\fF}$ implies the weak convergence 
of probability measures. We recall that a sequence 
$(P_{n})$ in $\cP(\Xi)$ converges weakly to $P$ iff
\[
\lim_{n\to\infty}\int_{\Xi}f(\xi)\,dP_{n}(\xi)=
\int_{\Xi}f(\xi)\,dP(\xi)
\]
holds for all bounded continuous functions $f:\Xi\to\R$.
Compared with classical probability metrics we consider 
here a much smaller class $\fF$ of functions, namely,
the collection of all integrands in \eqref{sopde}
\begin{equation}\label{integr}
\fF_{mi}=\{f(z,\cdot):z\in Z_{\rm ad}\}.
\end{equation} 
Following \cite{RaRo02} we call $d_{\fF_{mi}}$ the
problem-based or minimal information (m.i.) distance. 
Convergence with respect to such distances does not
imply weak convergence in general, but the question
arises whether it is implied by the weak convergence of 
probability measures. It is known that a positive answer 
depends on analytical properties of the class $\fF$. The 
following result is classical and due to \cite{Tops67}.
\begin{lem}\label{l1}
If $\fF$ is uniformly bounded and it holds that
\[
P(\{\xi\in\Xi:\fF\mbox{ is not equicontinuous at }
\xi\})=0,
\] 
then the set $\fF$ is a so-called $P$-uniformity class,
i.e., weak convergence of $(P_{n})$ to $P$ implies 
\[
\lim_{n\to\infty}d_{\fF}(P_{n},P)=0.
\]
\end{lem}
The choice \eqref{integr} of $\fF$ leads to the following 
result proved in \cite{HoRS20}. 
\begin{thm}\label{t1}
Under the standing assumptions and with the class $\fF_{mi}$ 
in \eqref{integr} we obtain the estimates 
\begin{eqnarray}\label{infstab}
|v(Q)-v(P)|&\leq& d_{\fF_{mi}}(P,Q)\\
\|z(Q) - z(P)\|_{H}&\leq& 2\sqrt{\frac{2}{\alpha}}
d_{\fF_{mi}}(P,Q)^{\frac12} 
\label{solstab}
\end{eqnarray}
for the optimal value $v(P)$ and solution $z(P)$ of
\eqref{sopde} if the original probability distribution 
$P$ is perturbed by any $Q\in\cP(\Xi)$.
\end{thm}
Next we collect some properties of the class $\fF_{mi}$ and 
of its elements under Lipschitz continuity assumptions 
on the coefficients of the linear elliptic PDE and of 
its right-hand sides implying that $\fF_{mi}$ is a
$P$-uniformity class (proved in \cite{HoRS20}).
\begin{thm}\label{t2}
Assume that all functions $b_{ij}(x,\cdot)$, 
$i,j=1,\ldots,m$, and $g(x,\cdot)$ are Lip\-schitz 
continuous on $\Xi$ uniformly with respect to $x\in D$, 
and let $g\in L_{\infty}(\Xi,P;H)$.
Then the family $\fF_{mi}=\{f(z,\cdot):z\in Z_{\rm ad}\}$ 
is uniformly bounded and Lipschitz continuous 
on $\Xi$ (with a constant not depending on $z$). 
In particular, $\fF_{mi}$ is a $P$-uniformity class. \\
\end{thm}
We close this section by extending the H\"older stability 
result \eqref{solstab} in Theorem \ref{t1} to Lipschitz
stability with respect to a pseudo-metric of the type
\eqref{zetadist}, but with a class $\fF_{di}$ of 
functions different from $\fF_{mi}$. For deriving the
Lipschitz stability result we do not make use of
classical work like, e.g., \cite{Alt91,DoHa93}, but
exploit the fact that \eqref{sopde} is formulated as
an optimization problem with fixed constraint set. Our
methodology exploits the quadratic growth condition 
of $F_{P}$ (see Proposition \ref{p1}) and partly 
parallels that of \cite[Lemma 2.1]{Shap92}.
\begin{thm}\label{t3}
Under the standing assumptions the Lipschitz-type estimate
\begin{equation}\label{Lipest}
\|z(Q)-z(P)\|_{H}\leq\frac{8}{\alpha}d_{\fF_{di}}(P,Q)
\end{equation}
holds for all $P,Q\in\cP(\Xi)$, where $\fF_{di}$ 
denotes the following function class on $\Xi$
\begin{equation}\label{derivintegr}
\fF_{di}=\left\{\langle A(\cdot)^{-1}(z+g(\cdot))-
\tilde{u},A(\cdot)^{-1}h\rangle_{H}+\alpha\langle 
z,h\rangle_{H}:z\in Z_{\rm ad},\|h\|_{H}\le 1\right\}.
\end{equation}
\end{thm}
\begin{pf}
Let $P,Q\in\cP(\Xi)$ and $z(P),z(Q)\in Z_{\rm ad}$ the
corresponding solutions to \eqref{sopde}. From Proposition
\ref{p1} we know that $F_{P}$ has quadratic growth 
around $z(P)$, i.e., 
\begin{eqnarray*}
\frac{\alpha}{8}\|z(Q)-z(P)\|_{H}^{2}&\leq& F_{P}(z(Q))-
F_{P}(z(P))\\
&\leq&(F_{P}(z(Q))-F_{Q}(z(Q)))-(F_{P}(z(P))-F_{Q}(z(P))),
\end{eqnarray*}
where we added $F_{Q}(z(P))-F_{Q}(z(Q))\geq 0$ to the
right-hand side. Now we consider the function
$h:[0,1]\to\R$ given by $h(t)=(F_{P}-F_{Q})(z(P)+
t(z(Q)-z(P)))$, $t\in[0,1]$. Due to Proposition \ref{p1}
$F_{P}$ and $F_{Q}$ are G\^ateaux differentiable on $H$ 
with G\^ateaux derivatives $F'_{P}$ and $F'_{Q}$.
Hence, $h$ is differentiable on $[0,1]$ and it holds that
$|h(1)-h(0)|\le\sup_{t\in[0,1]}|h'(t)|$. This implies
\begin{eqnarray*}
(F_{P}-F_{Q})(z(Q))-(F_{P}-F_{Q})(z(P))\!\!\!&\le&\!\!
\sup_{z\in Z_{\rm ad}}|(F'_{P}-F'_{Q})(z)(z(Q)-z(P))|\\
&\le&\!\!\!\sup_{z\in Z_{\rm ad}}\|(F'_{P}-F'_{Q})(z)\|
\|z(Q)-z(P)\|_{H}\\
\end{eqnarray*}
We obtain after dividing by $\|z(Q)-z(P)\|_{H}$ 
\begin{eqnarray*}
\frac{\alpha}{8}\|z(Q)-z(P)\|_{H}\! 
&\leq&\!\sup_{z\in Z_{\rm ad}}\|(F'_{P}-F'_{Q})(z)\|\\
&\leq&\!\sup_{z\in Z_{\rm ad}}\sup_{\|h\|_{H}\le 1}
\left|\int_{\Xi}f'_{z}(z,\xi)(h)\,d(P-Q)(\xi)\right|,
\end{eqnarray*}
where $f'_{z}(z,\xi)$ denotes the partial G\^ateaux
derivative of $f$ with respect to the first variable.
A straightforward evaluation shows that
\[
f'_{z}(z,\xi)(h)=\langle A(\xi)^{-1}(z+g(\xi))-\tilde{u},
A(\xi)^{-1}h\rangle_{H}+\alpha\langle z,h\rangle_{H}
\]
holds for all $z\in Z_{\rm ad}$, $\|h\|_{H}\le 1$,
$\xi\in\Xi$. This completes the proof.\hfill$\Box$
\end{pf}
\begin{rmk}\label{r1}
An inspection of the proof of Theorem \ref{t2} (see
\cite[Section 6]{HoRS20}) reveals that the class $\fF_{di}$
of (partial) derivatives of integrands in $\fF_{mi}$ is 
also a $P$-uniformity class under the assumptions of 
Theorem \ref{t2}.
\end{rmk}

\section{Monte Carlo approximations}
\label{sec:mc}

Let $\xi_1,\xi_2,\ldots,\xi_n,\ldots$ be independent 
identically distributed (iid) $\Xi$-valued random 
variables on some complete probability space 
$(\Omega,\cF,\W)$ having the common distribution $P$, 
i.e., $P=\W\xi_1^{-1}$. 
We consider the empirical measures
\begin{equation}\label{empmeas}
P_n(\cdot)=\frac{1}{n}\sum_{i=1}^n\delta_{\xi_i(\cdot)}
\qquad(n\in\N),
\end{equation}
where $\delta_{\xi}$ denotes the Dirac measure, which
places mass $1$ at $\xi\in\Xi$ and mass $0$ elsewhere.
Based on empirical measures we study the sequence of
empirical or Monte Carlo approximations of the stochastic
program \eqref{sopde} with sample size $n$, i.e.,
\begin{equation}\label{MCappr}
\min\left\{\int_{\Xi}f(z,\xi)\,dP_{n}(\cdot)(\xi)=
\frac{1}{n}\sum_{i=1}^nf(z,\xi_i(\cdot)):
z\in Z_{\rm ad}\right\}.
\end{equation}
The optimal value $v(P_{n}(\cdot))$ of \eqref{MCappr} is
a real random variable and the solution $z(P_{n}(\cdot))$ 
an $H$-valued random element 
(see \cite[Lemma III.39]{CaVa77}).\\
Qualitative and quantitative results on the asymptotic
behavior of optimal values and solutions to (\ref{MCappr})
are known in finite-dimensional settings 
(see \cite{DuWe88}, and the surveys \cite{Shap03} 
and \cite{Pflu03}). Since the sequence $(P_{n}(\cdot))$ 
of empirical measures converges weakly to $P$ $\W$-almost 
surely, one obtains the following corollary by combining 
Lemma \ref{l1} and Theorems \ref{t1} and \ref{t2}. 
\begin{cor}
The sequences $(v(P_{n}(\cdot)))$ and $(z(P_{n}(\cdot)))$ 
of empirical optimal values and solutions converge 
$\W$-almost surely to the true optimal values and 
solutions $v(P)$ and $z(P)$, respectively.
\end{cor}
In this section we are mainly interested in quantitative
results on the asymptotic behavior of $v(P_n(\cdot))$ 
and $z(P_n(\cdot))$. This is closely related to uniform 
convergence properties of the empirical process
\begin{equation}\label{empproc}
\left\{\G_n(\cdot)f:=\sqrt{n}(P_n(\cdot)-P)f=
\frac{1}{\sqrt{n}}\sum\limits_{i=1}^n(f(\xi_i(\cdot))-P f)
\right\}_{f\in\fF}
\end{equation}
indexed by a class $\fF$ of real-valued measurable 
functions on $\Xi$ and, hence, to quantitative estimates of
\begin{equation}\label{supemp}
\|\G_n(\cdot)\|_{\fF}=\sup_{f\in\fF}|\G_{n}(\cdot)f|=
\sqrt{n}\,d_{\fF}(P_n(\cdot),P)
=\sqrt{n}\,\sup_{f\in\fF}|P_n(\cdot)f - P f|.
\end{equation} 
Here, we set $P f=\int_{\Xi} f(\xi)\mathrm{d}P(\xi)$ for 
any probability distribution $P$ and any $f\in\fF$. 
Since the supremum in \eqref{supemp} is taken with respect
to an uncountable set $\fF$, it is not necessarily 
measurable with respect to $\cF$. In our applications
to the classes $\fF_{mi}$ and $\fF_{di}$, however, the
supremum \eqref{supemp} can be taken with respect
to some countable subset of both classes since the
Hilbert space $H$ is separable and all functions $f$
are continuous. Hence, the supremum is measurable.

There exist two main approaches to derive quantitative
information on the asymptotic behavior of empirical
processes. The first consists in the use of concentration
inequalities (pioneered in \cite{Tala95} and presented 
in some detail in \cite{BoLM13}) with applications 
to bounding \eqref{supemp} in probability. The second 
relies on the notion of Donsker classes of functions with 
applications to limit theorems. In this paper we study
the second approach.

A collection $\fF$ of measurable functions on $\Xi$ 
is called {\em $P$-Donsker} if the empirical
process \eqref{empproc} converges in distribution to a 
tight random variable $\G$ in the space 
$\ell^{\infty}(\fF)$, where the limit process 
$\G=\{\G f:f\in\fF\}$ is a zero-mean Gaussian process 
with the covariance function
\[
E_{P}[\G f_1 \G f_2]=P[(f_1-Pf_1)(f_2-Pf_2)]\quad
(f_1,f_2\in\fF).
\]
The limit $\G$ is sometimes called a $P$-Brownian 
bridge process in $\ell^{\infty}(\fF)$.
\begin{rmk}\label{r2}
We will prove that $\fF=\fF_{mi}$ and $\fF=\fF_{di}$ are 
$P$-Donsker classes by showing that 
$\sqrt{n}\E[d_{\fF}(P_n(\cdot),P)]$ is bounded 
(see Proposition \ref{p2}). From this we deduce the 
following mean convergence rates 
\begin{equation}\label{MCrate}
\E[d_{\fF}(P_n(\cdot),P)]=O(n^{-\frac12}).
\end{equation}
Together with Theorems \ref{t1} and \ref{t3} this then 
leads to best possible mean convergence rates of Monte 
Carlo estimates for optimal values and solutions.
\end{rmk}
Whether $\fF$ satisfies the $P$-Donsker class property, 
depends on its size measured in terms of so-called 
bracketing or metric entropy numbers. To introduce these 
concepts, let $\fF$ be a subset of the linear normed 
space $L_p(\Xi,P)$ (for some $p\ge 1$) (of equivalence
classes) of measurable functions endowed with the norm
\[
\|f\|_{P,\,p}=(P |f|^p)^{\frac1p}=\left(\int_{\Xi}
|f(\xi)|^{p}\,dP(\xi)\right)^{\frac1p}. 
\]
Given a pair of functions $l,u\in L_p(\Xi,P)$, $l\le u$, 
a bracket $[l,u]$ is defined by $[l,u]=\{f\in L_p(\Xi,P):
l\le f\le u\}$. Given $\varepsilon>0$ the bracketing number 
$N_{[\,]}(\varepsilon,\fF,\|\cdot\|_{P,p})$ is the minimal 
number of brackets with $\|l-u\|_{P,p}<\varepsilon$ 
needed to cover $\fF$. The metric entropy number with
bracketing of $\fF$ is defined by
\[
H_{[\,]}(\varepsilon,\fF,\|\cdot\|_{P,p})=\log{
N_{[\,]}(\varepsilon,\fF,\|\cdot\|_{P,p})}.
\]
Both numbers are finite if $\fF$ is a totally bounded 
subset of $L_p(\Xi,P)$. A powerful result on empirical 
processes is the following 
(see \cite[Thm. A.2]{vdVa96}).
\begin{prop}\label{p2}
There exists a universal constant $C>0$ such that for
any class $\fF$ of measurable functions with envelope 
function $\hat{F}$ (i.e., $|f|\le\hat{F}$ for every 
$f\in\fF$) belonging to $L_2(\Xi,P)$ the estimate
\begin{equation}\label{Donskbound}
\E[\|\G_n\|_{\fF}]\le C\int_{0}^{1}
\sqrt{1+H_{[\,]}(\varepsilon\|\hat{F}\|_{P,2},\fF,
\|\cdot\|_{P,2})}\,d\varepsilon\,\|\hat{F}\|_{P,2}
\end{equation}
holds. If the integral in \eqref{Donskbound} is 
finite, then the class $\fF$ is $P$-Donsker.
\end{prop}
Note that the integral in \eqref{Donskbound} can only
be finite if $H_{[\,]}(\varepsilon,\fF,\|\cdot\|_{P,2})$
grows at most like $\varepsilon^{-\beta}$ with $0<\beta<2$
for $\varepsilon\to +0$.\\
Next we discuss the assumption of finiteness of the
integral in \eqref{Donskbound} in case that $\fF$ is
a bounded subset of classical linear normed spaces of 
smooth functions.

\begin{exam}\label{e1}
Let $\Xi\subset\R^d$ be convex, bounded with the property
$\Xi\subseteq{\rm cl}\;{\rm int}\,\Xi$, $k\in\N_0$ and 
$r\in[0,1]$. We consider the linear space $C^{k,r}(\Xi)$ 
of real functions on $\Xi$ having partial derivatives up 
to order $k$ such that all $k$th order derivatives are 
H\"older continuous with exponent $r$. 
Next we use the notation ${\bf i}=(i_1,\ldots,i_d)$ 
with $i_j\in\N_0$, $j=1,\ldots,d$, and 
$|{\bf i}|=\sum_{j=1}^{d}i_j$. Further, $D^{\bf i}f$ denotes
\[
D^{\bf i}f=\frac{\partial^{|i|}f}{\partial\xi_{1}^{i_1}
\cdots\partial\xi_{d}^{i_d}}\qquad(f\in C^{k,r}(\Xi),
|{\bf i}|\le k).
\]
If the spaces are endowed with the norms
\begin{eqnarray*}
\|f\|_{k,0}&=&\max_{|{\bf i}|\le k}\sup_{\xi}
|D^{{\bf i}}f(\xi)|\\
\|f\|_{k,r}&=&\max_{|{\bf i}|\le k}\sup_{\xi}
|D^{\bf i}f(\xi)|+\max_{|{\bf i}|=k}
\sup_{\xi\neq\tilde{\xi}}\frac{|D^{\bf i}f(\xi)-
D^{\bf i}(\tilde{\xi})|}{\|\xi-\tilde{\xi}\|^r}\quad(r>0),
\end{eqnarray*}
where the suprema are taken over all $\xi,\tilde{\xi}$
in the interior of $\Xi$, they become Banach spaces
The metric entropy with bracketing of balls 
$\B_{k,r}(\rho)$ around the origin with radius $\rho$ 
in $C^{k,r}(\Xi)$ is computed in \cite{KoTi61} with
respect to the uniform norm 
$\|\cdot\|_{0,0}=\|\cdot\|_{\infty}$. 
The authors show that there exists a constant 
$K > 0$ depending only on $d$, $k$, $r$, $\rho$ and the 
diameter of $\Xi$ such that we have for every 
$\varepsilon>0$
\begin{equation}\label{metrentr}
H_{[\,]}(\varepsilon\rho,\B_{k,r}(\rho),
\|\cdot\|_{P,2})\le K\varepsilon^{-\frac{d}{k+r}},
\end{equation}
where the result from \cite{KoTi61} was adapted to the
norm in $L_{2}(\Xi,P)$ (see also
\cite[Section 2.7.1]{vdVWe96}).
Hence, Proposition \ref{p2} can be utilized to show
that bounded subsets of $C^{k,r}(\Xi)$ are
$P$-Donsker if $\frac{d}{2}<k+r$. 
For the situation studied in Theorem \ref{t2} with 
$k=0$ and $r=1$ this means that bounded subsets
of $C^{0,1}(\Xi)$ are $P$-Donsker only for $d=1$.
Without imposing stronger smoothness conditions on 
the coefficients in the bilinear form $a$ (see
\eqref{bilinform}) and the right-hand side compared 
to Theorem \ref{t2}, the integral in \eqref{Donskbound} 
will not be finite. Hence, one cannot use Proposition
\ref{p2} to conclude that $\fF$ is $P$-Donsker.
\end{exam}
\begin{rmk}\label{r3}
Indeed a convergence rate for the sequence 
$(\E[d_{\fF}(P_n(\cdot),P)])$ as in \eqref{MCrate}
cannot be achieved if $\fF$ is the unit ball in 
$C^{0,1}(\Xi)$ for $d>1$. Then $d_{\fF}$ coincides
with the Wasserstein metric $W_1$ and the Fortet-Mourier 
metric $\zeta_1$ of order $1$ (see also
\cite[Section 4]{HoRS20}). Namely, it is shown in
\cite{DeSS13,FoGu15} that the Wasserstein distance $W_p$
of $P$ and $P_n(\cdot)$ has the mean convergence rate
\begin{equation}\label{Wrate} 
\E[W_{p}(P_n(\cdot),P)]=O(n^{-\frac{1}{d}})
\end{equation}
if $d>2$, $p\ge 1$ and sufficiently high moments of $P$ 
exist. This rate carries over to the mean convergence
rate of Fortet-Mourier metrics $\zeta_p$ and of the
bounded Lipschitz metric $\beta$ which represents a
lower bound of $\zeta_1$ (see also \cite{Dudl69} for 
the mean convergence rate of the sequence 
$(\E[\beta(P_n(\cdot),P)])$).
\end{rmk}
Next we derive conditions implying that the functions
in $\fF_{mi}$ and $\fF_{di}$, respectively, are
sufficiently smooth. A similar result on differentiability
of solutions to random PDEs is proved in  
\cite[Section 4]{CoDS10} in a different way.
\begin{thm}\label{t4}
Let $\Xi\subset\R^d$ be a bounded, convex set having  
the property $\Xi\subseteq{\rm cl}\;{\rm int}\,\Xi$ and 
let $k\in\N$. Let the assumptions of Theorem \ref{t2} 
be satisfied and assume that, for all $u,v\in V$, the 
functions $\langle A(\cdot)u,v\rangle$ and 
$\langle g(\cdot),v\rangle$ belong to $C^{k,0}(\Xi)$.
Then both classes $\fF_{mi}$ and $\fF_{di}$ are subsets
of $C^{k,0}(\Xi)$.
\end{thm}
\begin{pf}
The integrands $f$ belonging to $\fF_{mi}$ are of the form
(see \eqref{sopde1})
\[
f(z,\xi)=\frac{1}{2}\left\|u(\xi)-\widetilde{u}
\right\|_{H}^{2}+\frac{\alpha}{2}\|z\|_{H}^{2}\,,
\]
where $u(\xi)=A(\xi)^{-1}(z+g(\xi))$, $\xi\in\Xi$. 
We begin by showing that, for $w,y\in V^{\star}$ the 
mapping $\xi\mapsto\langle y,A(\xi)^{-1}w\rangle$ has
first order partial derivatives at any 
$\xi\in{\rm int}\,\Xi$.
We fix $\xi\in{\rm int}\,\Xi$, $j\in\{1,\ldots,d\}$, and
a canonical basis vector $e_j\in\R^{d}$. Then
$\xi+he_j\in{\rm int}\,\Xi$ for sufficiently small 
$|h|>0$ and we obtain for $w, y\in V^{\star}$
\begin{eqnarray*}
\frac{1}{h}\langle y,(A(\xi+he_j)^{-1}-
A(\xi)^{-1})w\rangle\!&=&
\!\frac{1}{h}\langle y,A(\xi)^{-1}(A(\xi)-A(\xi+he_j))
u(h)\rangle\\
&=&\!\frac{1}{h}\langle\triangle_{A}^{j}(\xi;h)u,v\rangle-\\
&&\!\frac{1}{h}\langle(\triangle_{A}^{j}(\xi;h))^{\star}
v,u - u(h)\rangle,
\end{eqnarray*}
where $\triangle_{A}^{j}(\xi;h)=A(\xi+he_j)-A(\xi)$,
$u(h)=A(\xi+he_j)^{-1}w$,
$u=A(\xi)^{-1}w$, $v=(A(\xi)^{-1})^{\star}y\in V$ and 
$(A(\xi)^{-1})^{\star}$ denotes the adjoint mapping to
$A(\xi)^{-1}$. While the first summand on the right-hand
side converges for $h\to 0$ to the partial derivative 
$\frac{\partial}{\partial\xi_j}$ of $\langle A(\cdot)u,
v\rangle$ at $\xi$, the second converges to zero as 
$(u-u(h))$ converges to zero.
Hence, the partial derivative $\frac{\partial}
{\partial\xi_j}$ of $\langle y,A(\cdot)^{-1}w\rangle$
exists at $\xi$ and it holds
\[
\frac{\partial}{\partial\xi_j}\langle y,A(\xi)^{-1}w\rangle
=\frac{\partial}{\partial\xi_j}\langle A(\cdot)u,v
\rangle\quad(\mbox{at }\xi).
\]
This identity also shows $\langle y,A(\xi)^{-1}w\rangle$ 
is continuously differentiable. The differentiability of 
$\langle y,A(\cdot)^{-1}g(\cdot)\rangle$ follows in a
straightforward way via the product rule. Hence, we
conclude that the partial derivative
$\frac{\partial}{\partial\xi_j}f(z,\cdot)$ exists for any
$z\in Z_{\rm ad}$. By the same reasoning we can inductively
derive the existence of higher order mixed partial
derivatives $D^{\bf i}f(z,\cdot)$ at $\xi$ for 
$|{\bf i}|\le k$ and any $z\in Z_{\rm ad}$. We conclude 
that both classes $\fF_{mi}$ and $\fF_{di}$ are subsets 
of $C^{k,0}(\Xi)$.
\hfill$\Box$
\end{pf}

\begin{rmk}\label{r4}
According to the definition of the mapping $A(\xi):
V\to V^{\star}$ we have
\begin{equation}\label{diff}
\langle A(\xi)u,v\rangle=\sum_{i,j=1}^{m}\int_{D}
b_{ij}(x,\xi)\frac{\partial u(x)}{\partial x_{i}}
\frac{\partial v(x)}{\partial x_{j}}\,dx
\end{equation}
for all pairs $(u,v)\in V$. Due to the uniform ellipticity
condition \eqref{unifellipt} we know that all functions
$b_{ij}$ are essentially bounded. If we assume that all
functions $b_{ij}(x,\cdot):\Xi\to\R$, $x\in D$, have 
continuous mixed partial derivatives up to order $k$ 
which are in addition all measurable and essentially 
bounded on $D\times\Xi$, one obtains mixed partial 
derivatives of $\langle A(\cdot)u,v\rangle$ by
differentiating equation \eqref{diff}.\\
The same is true for $\langle g(\cdot),v\rangle$ if 
the functions $g(x,\cdot)$, $x\in D$, have continuous 
mixed partial derivatives up to order $k$ which are all 
measurable and essentially bounded on $D\times\Xi$.
\end{rmk}
In order to make use of Example \ref{e1} we present
conditions implying that both classes $\fF_{mi}$ and
$\fF_{di}$ are bounded subsets of $C^{k,0}(\Xi)$.

\begin{thm}\label{t5}
Let $\Xi\subset\R^d$ be a bounded, convex set having  
the property $\Xi\subseteq{\rm cl}\;{\rm int}\,\Xi$ 
and let $k\in\N$ be such that $d<2k$.
Assume that all functions $b_{ij}(x,\cdot):\Xi\to\R$,
$i,j=1,\ldots,m$, and $g(x,\cdot):\Xi\to\R$, $x\in D$,
have continuous partial derivatives up to order $k$ 
which are all measurable and essentially bounded 
on $D\times\Xi$.
Then the classes $\fF_{mi}$ and $\fF_{di}$ are 
$P$-Donsker, it holds that
\begin{eqnarray}\label{rateinf}
\E[|v(P_{n}(\cdot))-v(P)|]&=&O(n^{-\frac12})\\
\E[\|z(P_{n}(\cdot))-z(P)\|_{H}]&=&O(n^{-\frac12})
\label{ratesol}
\end{eqnarray}
and the sequence $(\sqrt{n}(v(P_n(\cdot))-v(P)))$
converges in distribution to some real random variable
$\zeta$, where $v(P)$ and $z(P)$ are the optimal value 
and solution of \eqref{sopde}, and $v(P_n(\cdot))$ and
$z(P_n(\cdot))$ are the optimal value and solution of \eqref{MCappr}, respectively.
\end{thm}
\begin{pf}
Our assumptions together with Theorem \ref{t4} imply that 
both classes $\fF_{mi}$ and $\fF_{di}$ represent bounded
subsets of the Banach space $C^{k,0}(\Xi)$. Hence,
according to Example \ref{e1} the metric entropy with
bracketing of any of the two classes satisfies
\[
H_{[\,]}(\varepsilon\rho,\fF,\|\cdot\|_{P,2})
\le K\varepsilon^{-\frac{d}{k}}
\]
for some constant $K>0$. Since $\Xi$ is bounded, the
estimate \eqref{Donskbound} in Proposition \ref{p2} implies
\[
\E[\sqrt{n}d_{\fF}(P_n(\cdot),P)]\le\hat{C}\int_{0}^{1}
\varepsilon^{-\frac{d}{k}}d\varepsilon
\]
for some $\hat{C}>0$. Since $\frac{d}{k}<2$, the 
right-hand side is bounded and we have that 
\[
\E[d_{\fF}(P_n(\cdot),P)]=O(n^{-\frac12})
\]
holds for $\fF=\fF_{mi}$ and $\fF=\fF_{di}$. Hence, we
obtain \eqref{rateinf} from Theorem \ref{t1} and
\eqref{ratesol} from Theorem \ref{t3}. Furthermore,
we conclude for $\fF=\fF_{\rm mi}$ from Proposition 
\ref{p2} that the empirical process 
$\{\G_n(\cdot)f=\sqrt{n}(P_n(\cdot)-P)f\}_{f\in\fF}$
converges in distribution to a tight random variable 
$\{\G f\}_{f\in\fF}$ on $(\Omega,\cF,\W)$ with  
values in the space $\ell^{\infty}(\fF)$.\\ 
Due to the structure \eqref{integr} of $\fF$, we may
also write $\{Pf\}_{f\in\fF}=\{Pf(z,\cdot):
z\in Z_{\rm ad}\}$ and $\{\G f\}_{f\in\fF}=
\{\G f(z,\cdot)\}_{z\in Z_{\rm ad}}$.
This means that $Pf$ may be considered as element of 
the space $\ell^{\infty}(Z_{\rm ad})$ of bounded real-valued 
functions on $Z_{\rm ad}$. Correspondingly, $\G f$ may be
viewed as random variable in $\ell^{\infty}(Z_{\rm ad})$.
It remains to utilize the functional delta theorem 
(see \cite{Roem06}) for the infimal mapping 
\[
\Phi:\ell^{\infty}(Z_{\rm ad})\to\R,\quad\Phi(h)=
\inf_{z\in Z_{\rm ad}}h(z).
\]
The mapping $\Phi$ is finite, concave, hence, directionally
differentiable on $\ell^{\infty}(Z_{\rm ad})$. 
In addition, $\Phi$ is Lipschitz continuous (with modulus 
$1$) and, hence, Hadamard directionally differentiable 
(see \cite{Shap90}). The Hadamard directional derivative 
at $h_0\in\ell^{\infty}(Z_{\rm ad})$ is of the form
(see \cite{Lach06} and the survey \cite{Roem06})
\begin{equation}\label{Hadamderiv}
\Phi'_{h_0}(h)=\lim_{\varepsilon\downarrow 0}\inf\{h(z):
z\in Z_{\rm ad},h_0(z)\le\Phi(h_0)+\varepsilon\}
\quad(h\in\ell^{\infty}(Z_{\rm ad})).
\end{equation}
Then the functional delta theorem \cite[Theorem 1]{Roem06}
implies that 
\begin{equation}\label{funcdelta} 
\sqrt{n}(\Phi(P_n(\cdot)f)-\Phi(Pf)))\stackrel{d}
\longrightarrow\Phi'_{Pf}(\G f)=\zeta,
\end{equation}
where $\stackrel{d}\longrightarrow$ denotes convergence in
distribution of real random variables. Since the 
Hadamard directional derivative $\Phi'_{Pf}(\cdot)$ is
continuous (see \cite{Shap90}), $\zeta$ is a real
random variable on $(\Omega,\cF,\W)$. This completes 
the proof.\hfill$\Box$
\end{pf}
To establish an extension of Theorem \ref{t5} to 
$\Xi=\R^d$ let $\R^d=\bigcup_{j=1}^{\infty}\Xi_j$ be a 
partition of $\R^d$, where each set $\Xi_j$ is bounded, 
convex and has the property 
$\Xi_j\subseteq{\rm cl}\;{\rm int}\,\Xi_j$, $j\in\N$.
The idea is to apply Theorem \ref{t5} on each subset
$\Xi_j$ of $\R^d$ and then to apply the argument in
\cite[Theorem 1.1]{vdVa96} (see also
\cite[Corollary 2.1]{vdVa96}).
\begin{cor}\label{c2}
Let $k\in\N$ be such that $d<2k$. Assume that all 
functions $b_{ij}(x,\cdot):\R^d\to\R$, $i,j=1,\ldots,m$, 
and $g(x,\cdot):\R^d\to\R$, $x\in D$, have continuous
partial derivatives up to order $k$ which are all 
measurable on $D\times\Xi$.
Moreover, assume that for each $j\in\N$ the restrictions 
to $\Xi_j$ of all functions in both classes $\fF_{mi}$
and $\fF_{di}$ belong to the ball $\B_{k,0}(\rho_j)$ 
in $C^{k,0}(\Xi_j)$ and that the probability measure 
$P$ satisfies
\begin{equation}\label{tailcond}
\sum_{j=1}^{\infty}\rho_jP(\Xi_j)^{\frac12}<\infty.
\end{equation}
Then both classes $\fF_{mi}$ and $\fF_{di}$ are 
$P$-Donsker and \eqref{rateinf}, \eqref{ratesol} and
the central limit theorem for optimal values remain true.
\end{cor}
We note that condition \eqref{tailcond} represents a quite
implicit link between the growth of derivatives of the
functions in both classes with the tail behaviour of $P$.
\begin{rmk}\label{r5}
Note that Theorem  \ref{t5} allows to derive
asymptotically consistent confidence intervals for
optimal values by using resampling techniques 
such as bootstrapping \cite{GiZi90} or subsampling 
\cite{PoRW99}. Since we not know that the Hadamard
directional derivative is linear in the direction,
the classical bootstrap cannot be used. However, 
a variant called extended bootstrap in \cite{EiRo07}
and also subsampling can be used. The subsampling 
method is more generally applicable than the 
bootstrap, because only a basic limit theorem like 
that in Theorem \ref{t5} is required.
\end{rmk}

\section{Subsampling}
\label{sec:subsampl}

The subsampling method \cite{PoRW99} is based on 
sampling and resampling, but resampling is performed 
repeatedly without replacement and with a lower sample 
size $b=b(n)\in\N$, $b\ll n$. For some sufficiently 
large $n$, let $\xi_{1},\ldots,\xi_{n}$ be an iid sample
from $P$. Let $P_n$ be the empirical measure and $v(P_n)$
the corresponding optimal value of \eqref{MCappr}.
Based on the samples $\xi_{n_1},...,\xi_{n_b}$ drawn 
from $\{1,...,n\}$ with cardinality $b$, we
consider the corresponding empirical measure
\[
P^*(n_1,...,n_b)=\frac{1}{b}\sum_{i=1}^{b}
\delta_{\xi_{n_i}}
\]
and the optimal value $v(P^*(n_1,...,n_b))$.
The subsampling method estimates the limit distribution 
of $\zeta=\Psi'_{Pf}(\G f)$ (see \eqref{funcdelta}) based 
on both optimal values. It is justified by the limit theorem 
\cite[Theorem 2.1]{PoRo94} which reads in our framework
\begin{equation}\label{subsam1}
\left(\begin{array}{c} n \\ b \end{array} \right)^{-1}
\sum_{1\le n_1<...<n_b\le n}
\delta_{ \big\{
\sqrt{b} \big( v(P^*(n_1,...,n_b)) - v(P_n) \big)
\big\} }
\stackrel{d}\longrightarrow\zeta
\end{equation}
for $b,n\rightarrow\infty$ and $b/n\rightarrow 0$. 
The number of summands in \eqref{subsam1} becomes 
extremely large as $n$ and $b$ grow. However, the 
result remains valid if a number $m=m(n)$ is chosen
and the sum over all possible subsets is replaced by
the sum over $m$ randomly chosen subsets of 
$\{1,...,n\}$ of cardinality $b$: 
Let $N_j^{n,b}\subset\{1,...,n\}$ be randomly chosen 
with cardinality $\#N_j^{n,b}=b$ for $j=1,...,m$. 
Then with $P_n^*(N_j^{n,b})$ denoting the empirical 
measure based on $\{\xi_i : i\in N_j^{n,b}\}$ we have
\begin{equation}\label{subsam2}
\frac{1}{m}\sum_{j=1}^m
\delta_{\left\{\sqrt{b} 
\left( v(P_n^*(N_j^{n,b})) - v(P_n)\right)\right\} }
\stackrel{d}\longrightarrow\zeta
\end{equation}
for $b,n,m\rightarrow\infty$ and $b/n\rightarrow 0$ 
\cite[Corollary 2.1]{PoRo94}.

This suggests the following procedure to determine
confidence intervals for $v(P)$: Given $n,b,m\in\N$, 
$b<n$, sufficiently large and a sample 
$\xi_{1},\xi_{2},...,\xi_{n}$ from $P$. 
Compute $v(P_n)$. Resample from $P_n$ without
replacement with sample size $b<n$ to obtain
$\{\xi_i:i\in N_j^{n,b}\}$. Compute $v(P_n^*(N_j^{n,b}))$
and repeat this $m$ times. Let
\[
L_m(t)=\frac{1}{m}\sum_{j=1}^m
\one_{\left\{\sqrt{b} \left( v(P_n^*(N_j^{n,b})) - v(P_n)\right)\le t\right\} }
\quad(t\in\R).
\]
Choose $\alpha\in(0,1)$ and calculate the quantile
$\zeta_{1-\alpha,m}^{*}=\inf\{t:L_m(t)\ge 1-\alpha\}$ 
of $L_m$. Then we obtain for the asymptotic coverage 
probability of $v(P)$ 
\[
\lim_{n,m\to\infty}\W\left\{\sqrt{n}(v(P_{n})-v(P))\le
\zeta_{1-\alpha,m}^{*}\right\}\ge 1-\alpha.
\]

\section{Discussion and conclusions}
\label{sec:concl}

In this paper we studied Monte Carlo methods for solving
a stochastic optimization problem with linear quadratic
risk-neutral objective function, a linear elliptic PDE 
with random coefficients and convex control constraints.
Based on empirical process theory we were able to show 
that both optimal values and solutions converge in mean
with the best possible convergence rate $O(n^{-\frac12})$
if the coefficients of the PDE are sufficiently smooth.
The required degree of smoothness is related to the 
finite dimension of the random parameter. In addition, 
the optimal values satisfy a central limit result 
which enables the derivation of confidence intervals 
by resampling.

Our methodology is no longer successful if the 
optimization model \eqref{sopde} contains random convex 
control constraints that correspond to state constraints 
in the original stochastic optimization problem
\eqref{object}, \eqref{rellpde}. It also fails 
if the risk-neutral expectation in the objective
is replaced by some convex risk measure. Although
such risk measures preserve convexity, they typically
introduce nonsmoothness as, for example, in the case 
of so-called Conditional or Average Value-at-risk 
${\rm CVaR}$. In this case, problem \eqref{sopde} 
would be of the form
\begin{equation}\label{cvarsopde}
\min\Big\{{\rm CVaR}_{\kappa}(f(z,\cdot))=
\inf_{t\in\R}\Big\{t+\frac{1}{1-\kappa}\!\int_{\Xi}
\max\{0,f(z,\xi)-t\}\,dP(\xi)\Big\}:z\in Z_{\rm ad}
\Big\}
\end{equation}
for some $\kappa\in(0,1)$ and $f$ defined in 
\eqref{sopde1}. Hence, the corresponding 
minimal information distance is based on a class $\fF$ 
of functions that is no longer smooth as needed for our 
main result (Theorem \ref{t5}). The classical way for
reformulating \eqref{cvarsopde} into a smooth optimization
problem was suggested in \cite{RoUr02} and leads to
\[
\min\Big\{t+\frac{1}{1-\kappa}\int_{\Xi}y(\xi)
\,dP(\xi):y(\xi)\ge f(z,\xi)-t,\,y(\xi)\ge 0,
\,t\in\R_{+},\,z\in Z_{\rm ad}\Big\}
\]
and, thus, to an optimization model with random convex 
constraints. A possible way out consists in the approach
of smoothing CVaR as suggested in \cite{KoSu16,KoSu19}.

Finally, we mention two possible extensions of the 
results in this paper. The first extension consists in 
introducing a random mapping $B(\xi):H\to V^{\star}$ and
by replacing $z$ in \eqref{rellpde} by $B(\xi)z$. If
one requires that the function $\langle B(\cdot)z,v\rangle$
is sufficiently smooth on $\Xi$ for all $z\in Z_{\rm ad}$, 
$v\in V$, the function classes $\fF_{mi}$ and $\fF_{di}$ 
have to be modified, but the main results carry over.\\
A second extension concerns the finite dimensionality 
of $\Xi$ in Theorems \ref{t4} and \ref{t5}. In our earlier
paper \cite{HoRS20} and in Section \ref{sec:problem} the 
set $\Xi$ represents a metric space. Hence, the general
stability results (Theorems \ref{t1} and \ref{t3})
enable the use of probability measures on infinite 
dimensional spaces. For example, this allows to consider 
the Karhunen-Lo\`eve expansion of a centered stochastic 
process $\{\xi_x\}_{x\in D}$ with probability distribution
$P$ on $\Xi=L_2(D)$, finite second moments and continuous
covariance function $K(x,y)=\E[\xi_x\,\xi_y]$, $x,y\in D$,
which is of the form
\begin{equation}\label{KLexp}
\xi_x=\sum_{j=1}^{\infty}Z_je_j(x)\qquad(x\in D).
\end{equation}
Here, $(e_j)_{j\in\N}$ is an orthogonal system in $H=L_2(D)$
and $(Z_j)_{j\in\N}$ is a sequence of centered, 
uncorrelated real random variables (see \cite{Stei19} and
references therein). A truncated version of \eqref{KLexp}
with $d$ summands can then be used for Monte Carlo 
sampling and the truncation error be estimated by the 
distance $d_{\fF}$ or possible upper bounds.

\end{document}